\documentclass[12pt]{amsart}
\usepackage{amsmath,amssymb}
\usepackage{amsfonts}
\usepackage{amsthm}
\usepackage{latexsym}
\usepackage{graphicx}
\def\p{\partial}

\def\Hom{\mbox{Hom}}

\def\R{\mathbb{R}}

\def\vv<#1>{\langle#1\rangle}

\def\XXint#1#2{\setbox0=\hbox{$#1{#2}{\int}$}{#2}\kern-.5\wd0 }

\def\XXint#1#2#3{{\setbox0=\hbox{$#1{#2#3}{\int}$}
     \vcenter{\hbox{$#2#3$}}\kern-.5\wd0}}

\def\vv<#1>{\left\langle#1\right\rangle}
\def\pt{\frac{\p}{\p t}}
\def\pu{\frac{\p}{\p u}}
\def\ol{\overline}
\def\wh{\widehat}
\def\e{\epsilon}
\newtheorem{thm}{Theorem}[section]

\newtheorem{lem}{Lemma}[section]
\newtheorem{prop}{Proposition}[section]
\newtheorem{cor}{Corollary}[section]
\theoremstyle{definition}
\newtheorem{defn}{Definition}[section]
\theoremstyle{remark}

\numberwithin{equation}{section}
\def\dev{{\rm dev}}

\def\wt{\widetilde}
\begin{document}

\title{Extensions of isometric immersions and Cartan-Ambrose-Hicks theorem based on submanifolds}

\author{Chengjie Yu$^1$}
\address{Department of Mathematics, Shantou University, Shantou, Guangdong, 515063, China}
\email{cjyu@stu.edu.cn}

\thanks{$^1$Research partially supported by GDNSF with contract no. 2021A1515010264 and NNSF of China with contract no. 11571215.}

\renewcommand{\subjclassname}{%
  \textup{2010} Mathematics Subject Classification}
\date{}
\begin{abstract}
In this paper, we study the general extension problem for isometric immersions by establishing Cartan-Ambrose-Hicks theorems based on submanifolds. Our method also provides geometric constructions of such extensions.
\end{abstract}
\maketitle\markboth{Chengjie Yu }{Extensions of isometric immersions}
\section{Introduction}
Let $(S,g_S)$, $(M,g)$ and $(\wt M,\wt g)$ be Riemannian manifolds, and $\iota:S\to M$ and $\wt\iota:S\to \wt M$ be isometric immersions. In this paper, we consider the extension problem of isometric immersions: the existence of isometric immersions $f:M\to \wt M$ such that $f\circ \iota=\wt\iota$. When $\dim \wt M=\dim M$, we consider this as a Cartan-Ambrose-Hicks theorem \cite{Am,Ca,Hi} based on a submanifold and deal with it in a similar way to our previous work \cite{Yu} using developments of curves. When $\dim \wt M>\dim M$, we consider this as a positive codimensional Cartan-Ambrose-Hicks theorem based on submanifolds and deal with it in a similar way to our previous work \cite{Yu2} using developments of curves w.r.t. symmetric tensors.

We first recall the definition of developments of curves in \cite{Yu}.
\begin{defn}\label{def-dev}
Let $(M,g)$ be a Riemannian manifold, $p\in M$ and $v:[0,T]\to T_pM$ be a smooth curve in $T_pM$. The development of $v$ is a curve $\gamma:[0,T]\to M$ such that
$$\left\{\begin{array}{ll}\gamma(0)=p&\\
\gamma'(t)=P_0^t(\gamma)(v(t))& t\in [0,T].
\end{array}\right.$$
Here $P_{t_1}^{t_2}(\gamma)$ is the parallel displacement from $\gamma(t_1)$ to $\gamma(t_2)$ along $\gamma$.
\end{defn}

In the monograph \cite{KN} by Kobayashi and Normizu,  the curve $v(t)$ in Definition \ref{def-dev} is called the development of $\gamma(t)$ reversely. In \cite{Yu}, we give a detailed proof of the uniqueness and local existence for developments of curves. The global existence of developments of curves is equivalent to the completeness of the Riemannian metric (see \cite{KN}). We will denote the development of $v$ as $\dev(p,v)$ when it exists. Conversely, let  $\gamma:[0,1]\to M$ be a smooth curve with $\gamma(0)=p$, define
$$v_\gamma(t):=P_t^0(\gamma)(\gamma'(t)),$$
then $$\gamma=\dev(p,v_\gamma)$$ by Definition \ref{def-dev}.

Before stating the first main result, we give some explanations on notations and settings. Let $(M,g)$  and $(S,g_S)$ be Riemannian manifolds and $\iota:S\to M$ be an isometric immersion. For any $p\in S$, let
\begin{equation}\label{eq-orthogonal-decomp}
T_{\iota(p)}M=\iota_{*p}(T_pS)\oplus\left(\iota_{*p}(T_pS)\right)^\perp
\end{equation}
be the natural orthogonal decomposition. For a vector $v\in T_{\iota(p)}M$, we denote its orthogonal projection into $\iota_{*p}(T_pS)$ as $v^\top$ and denote its orthogonal projection into $\left(\iota_{*p}(T_pS)\right)^\perp$ as $v^\perp$.

 Let $\iota_*(TS)$ be the smooth vector bundle on $S$ with $$\iota_*(TS)_p=\iota_{*p}(T_pS),\ \forall p\in S,$$
 and $T^\perp S$ be the smooth vector bundle on $S$ with
 $$T_p^\perp S:=(T^\perp S)_p=\left(\iota_{*p}(T_pS)\right)^\perp,\ \forall p\in S.$$ Then, we have the orthogonal decomposition of smooth vector bundles:
$$\iota^*TM=\iota_*(TS)\oplus T^\perp S,$$
and the pull-back connection $\iota^*\nabla$ of the Levi-Civita connection of $\nabla$ of $g$ on $M$ induces canonical connections $\nabla^\top$ and $\nabla^\perp$ on the smooth vector bundles  $\iota_*(TS)$ and $T^\perp S$ respectively according to the orthogonal decomposition in the natural way:
$$\nabla_X^\top Y=((\iota^*\nabla)_XY)^\top,\ X\in \Gamma(TS),\ Y\in \Gamma(\iota_*(TS))$$
and
$$\nabla_X^\perp\xi=((\iota^*\nabla)_X\xi)^\perp,\ X\in \Gamma(TS),\ \xi\in \Gamma(T^\perp S).$$
A basic result in geometry of submanifolds says that the bundle isomorphism:$$\iota_*:TS\to \iota_*(TS)$$ preserves metrics and connections. So, for simplicity, we identify the two vector bundles by $\iota_*$. That is, a tangent vector of $S$ is also viewed as a tangent vector of $M$ by $\iota_*$. Moreover, a curve in $S$ is also viewed as a curve  in $M$ by composite it with $\iota$. We also denote the normal bundle $T^\perp S$ of the isometric immersion $\iota:S\to M$ as $N(\iota)$ and denote $T_p^\perp S$ as $N_p(\iota)$.

Let $\wt\iota:S\to (\wt M,\wt g)$ be another isometric immersion. Then, a homomorphism $\psi:N_p(\iota)\to N_p(\wt\iota)$ can be extended to a homomorphism
$\ol\psi:T_{\iota(p)}M\to T_{\wt\iota(p)}\wt M$ according to the orthogonal decomposition \eqref{eq-orthogonal-decomp} such that $\ol\psi|_{T_pS}$ is the identity map after the identifications of $\iota_{*p}(T_pS)$ and $\wt\iota_{*p}(T_pS)$ with $T_pS$.

We are now ready to state our first main result about the extension of isometric immersions in the case that $\dim M=\dim \wt M$.
\begin{thm}\label{thm-CAH-isometry}
Let $(M^n,g)$ and $(S^r,g_S)$ be two connected Riemannian manifolds and $\iota:S\to M$ be an  isometric immersion such that $\iota$ induces a surjective homomorphism between fundamental groups of $S$ and $M$. Let $(\wt M^n,\wt g)$ be a complete Riemannian manifold and $\wt \iota:S\to \wt M$ be an isometric immersion. Let  $\psi:N(\iota)\to N(\wt\iota)$ be a homomorphism of smooth vector bundles preserving metrics, second fundamental forms and normal connections. For any smooth  curve $\gamma:[0,1]\to M$ with $\gamma(0)=\iota(p)$ for $p\in S$, let $\wt v_\gamma=\ol\psi(v_\gamma)$ and $\wt \gamma=\dev(\wt\iota(p), \wt v_\gamma)$. Moreover, suppose that
\begin{equation}\label{eq-R-1}
R=\tau^{*}_\gamma \wt R
\end{equation}
for any smooth curve $\gamma:[0,1]\to M$ with $\gamma(0)=\iota(p)$ for $p\in S$ where $$\tau_\gamma=P_0^1(\wt \gamma)\circ \ol\psi\circ P_1^0(\gamma):T_{\gamma(1)}M\to T_{\wt\gamma(1)}\wt M,$$
and $R$ and $\wt R$ are the curvature tensors of $(M,g)$ and $(\wt M,\wt g)$ respectively. Then, the map $f(\gamma(1))=\wt\gamma(1)$ from $M$ to $\wt M$ is well defined and $f$ is the local isometry from $M$ to $\wt M$ with $\wt\iota=f\circ\iota$ and $f_{*}|_{N(\iota)}=\psi$.
\end{thm}

Note that when $M$ is simply connected, $\iota:S\to M$ automatically induces a surjective homomorphisms between fundamental groups of $S$ and $M$. So, one of the advantages in considering Cartan-Ambrose-Hicks theorem based on submanifolds is that one can relax the topological assumption of the source manifold from simply connectedness to relative simply connectedness.

Moreover, when taking $M=\wt M$ as space forms (simply connected Riemannian manifolds with constant sectional curvature), Theorem \ref{thm-CAH-isometry} will give us a rigidity result of isometric immersions with prescribed second fundamentals since the curvature assumption \eqref{eq-R-1} is clearly true in this case. This is a higher dimensional and higher codimensional extension of Bonnet's uniqueness theorem (see \cite{Bo,Chen}). We denote the space form of constant sectional curvature $K$ as $\mathbb M_K$.
\begin{cor}\label{cor-rigidity-isometry}
Let $(S^r,g_{S})$ be a connected Riemannian manifold, and $\iota,\wt\iota:S\to\mathbb M^n_K$ be two isometric immersions such that there is  a bundle map $\psi:N(\iota)\to N(\wt\iota)$ preserving metrics, second fundamental forms and normal connections. Then, there is an isometry $f:\mathbb M_K\to \mathbb M_K$ such that $\wt\iota=f\circ\iota$ and $f_*|_{N(\iota)}=\psi$.
\end{cor}

We next come to consider extensions of isometric immersions when $\dim \wt M>\dim M$. In this case, the second fundamental form of the isometric immersion of $M$ into $\wt M$ must be involved because of the Gauss-Codazzi-Ricci equations for submanifolds. We need to generalize the notion of developments of curves to the positive codimensional case which is the developments of curves w.r.t. symmetric tensors introduced in our previous work \cite{Yu2}.

\begin{defn}[Development of curve w.r.t. a symmetric tensor \cite{Yu2}]\label{def-g-dev}
Let $(\wt M^{n+s},\wt g)$ be a Riemannian manifold and $\wt p\in \wt M$. Let $T_{\wt p}\wt M=T^n\oplus N^s$ be an orthogonal decomposition,  $\wt h(t):[0,l]\to \Hom(T\odot T, N)$ and $\wt v:[0,l]\to T$ be smooth maps. A curve $\wt \gamma:[0,l]\to \wt M$ is called a development of $\wt v$ w.r.t. $\wt h$ if there exists a moving orthonormal frame $\left\{\wt E_A\Big|A=1,2,\cdots,n+s\right\}$ along $\wt\gamma $ satisfying the following equations:
\begin{equation}\label{eq-g-dev}
\left\{\begin{array}{ll}
\wt\gamma(0)=\wt p\\
\wt e_a:=\wt E_a(0)\in T&a=1,2,\cdots,n\\
\wt e_\alpha:=\wt E_\alpha(0)\in N&\alpha=n+1,n+2,\cdots,n+s\\
\wt \nabla_{\wt\gamma'(t)}\wt E_a=\displaystyle\sum_{\alpha=n+1}^{n+s}\vv<\wt h(t)(\wt v(t),\wt e_a),\wt e_\alpha>\wt E_\alpha& a=1,2,\cdots,n\\
\wt\nabla_{\wt\gamma'(t)}\wt E_\alpha=-\displaystyle \sum_{a=1}^n\vv<\wt h(t)(\wt v(t),\wt e_a),\wt e_\alpha>\wt E_a&\alpha=n+1,\cdots,n+s\\
\wt\gamma'(t)=\displaystyle\sum_{a=1}^n\vv<\wt v(t),\wt e_a>\wt E_a.\\
\end{array}\right.
\end{equation}
Here $T\odot T$ means the symmetric product of $T$. Moreover, define the map $D_{t_1}^{t_2}(\wt\gamma)$ as
 $$D_{t_1}^{t_2}(\wt \gamma):T_{\wt \gamma(t_1)}\wt M\to T_{\wt \gamma(t_2)}\wt M,\
 \sum_{A=1}^{n+s}c_A\wt E_A(t_1)\mapsto\sum_{A=1}^{n+s}c_A\wt E_A(t_2)$$
 for any $c_1,c_2\cdots,c_{n+s}\in \R$.
\end{defn}
In \cite{Yu2}, we prove the local existence and uniqueness of the developments of curves w.r.t.  symmetric tensors, and show its global existence when $\wt M$ is complete. We will then denote the curve $\wt \gamma$ in  Definition \ref{def-g-dev} as $\dev(\wt p,\wt v,\wt h)$ when it exists.

Before stating the our second main result about extension of isometric immersions in the positive codimensional case, we first fix some notations that will be used. Let $(M^n,g)$ be a Riemannian manifold,  $(V^s,\mathfrak{h},D)$ be a smooth vector bundle on $M$ equipped with the Riemannian metric $\mathfrak{h}$ and the compatible connection $D$, and $h\in \Gamma(\Hom(TM\odot TM,V))$. Let $\iota:(S^r,g_S)\to (M,g)$ be an isometric immersion. We will equip with the vector bundle $T^\perp S\oplus \iota^*V$ the natural direct-sum metric and the connection $\wh D$ defined as follows:
\begin{equation}\label{eq-hat-D-1}
\wh D_X\xi=\nabla^\perp_{X}\xi+h(X,\xi)\ (\forall\ X\in \Gamma(TS),\ \xi\in \Gamma(T^\perp S))
\end{equation}
and
\begin{equation}\label{eq-hat-D-2}
\wh D_X\eta=-(A_\eta(X))^\perp+D_X\eta\ (\forall\ X\in \Gamma(TS),\ \eta\in \Gamma(V)).
\end{equation}
Here $\nabla^\perp$ is the normal connection on $T^\perp S$ and for any $\eta\in V_p$,  $A_\eta:T_pM\to T_pM$  is defined by
\begin{equation*}
\vv<A_\eta(X),Y>_g=\vv<h(X,Y),\eta>_{\mathfrak{h}},\ \forall\ X,Y\in T_pM.
\end{equation*}
It is not hard to see that $\wh D$ is compatible with the natural direct-sum metric on $T^\perp S\oplus \iota^* V$.

Moreover, let $\wt\iota:S\to (\wt M^{n+s},\wt g)$ be another isometric immersion and $\psi:N_p(\iota)\oplus (\iota^*V)_p\to N_p(\wt\iota)$ be a linear map with $p\in S$.  Then, we denote $$\ol\psi:T_{\iota(p)}M\oplus V_{\iota(p)}=T_pS\oplus N_p(\iota)\oplus (\iota^*V)_p\to T_pS\oplus N(\wt\iota)=T_{\wt\iota(p)}\wt M$$ the orthogonal extension of $\psi$ such that $\ol\psi|_{T_pS}$ is the identity map.

We are now ready to state the second main result.
\begin{thm}\label{thm-CAH-immersion}
Let $(M^n,g)$ and $(S^r,g_{S})$ be two connected Riemannian manifolds, and $\iota:S\to M$ be an isometric immersion with second fundamental form $\sigma$ that induces a surjective homomorphisms bwtween fundamental groups of $S$ and $M$.  Let $(\wt M^{n+s},\wt g)$ be a complete Riemannian manifold and $\wt\iota :S\to \wt M$ be an isometric immersion with second fundamental form $\wt\sigma$. Let $(V^s,\mathfrak{h}, D)$ be a smooth vector bundle on $M$ equipped with the Riemannian metric $\mathfrak{h}$ and the compatible connection $D$, and $h\in \Gamma(\Hom(TM\odot TM, V))$.   Let $\psi:N(\iota)\oplus \iota^*V\to N(\wt\iota)$ be a bundle map that preserves metrics, connections and satisfies
\begin{equation}\label{eq-difference-1}
\psi^*\wt\sigma=\sigma+\iota^*h.
\end{equation}
For any smooth curve  $\gamma:[0,1]\to M$ with $\gamma(0)=\iota(p)$ for $p\in S$, let $\wt v_\gamma=\ol\psi(v_\gamma)$ and  $\wt\gamma=\dev(\wt\iota(p),\wt v_\gamma,\wt h_\gamma)$ where
$$\wt h_{\gamma}(t)=\left(\ol\psi^{-1}\right)^*\left(P_t^0(\gamma)h\right)\ \forall t\in [0,1].$$
Moreover, suppose that for any smooth curve $\gamma:[0,1]\to M$ with $\gamma(0)=\iota(p)$ for $p\in S$,\\
(1) for any $X,Y,Z,W\in T_{\gamma(1)}M$,
$$R(X,Y,Z,W)=\left(\tau_\gamma^*\wt R\right)(X,Y,Z,W)+\vv<h(X,W),h(Y,Z)>-\vv<h(X,Z),h(Y,W)>,$$
where $R$ and $\wt R$ are the curvature tensors of $M$ and $\wt M$ respectively;\\
(2) for any  $X,Y,Z\in T_{\gamma(1)}M$ and $\xi\in V_{\gamma(1)}$,
$$\vv<(D_{X}h)(Y,Z)-(D_Yh)(X,Z),\xi>=\left(\tau^*_\gamma\wt R\right)(Z,\xi,X,Y);$$
(3) for any $X,Y\in T_{\gamma(1)}M$ and $\xi,\eta\in V_{\gamma(1)}$,
$$R^V(\xi,\eta,X,Y)=\left(\tau_{\gamma}^*\wt R\right)(\xi,\eta,X,Y)+\vv<A_\xi(Y),A_\eta(X)>-\vv<A_\eta(Y),A_\xi(X)>,$$
 where $R^V$ is the curvature tensor of the vector bundle $V$. Here
    $$\tau_\gamma=D_0^1(\wt\gamma)\circ\ol\psi\circ P_1^0(\gamma):T_{\gamma(1)}M\oplus V_{\gamma(1)}\to T_{\wt\gamma(1)}\wt M.$$
Then, we have the following conclusions.
\begin{enumerate}
\item The map $f(\gamma(1))=\wt\gamma(1)$ from $M$ to $\wt M$ is well defined and is an isometric immersion such that $\wt\iota=f\circ\iota$ and $f_*|_{N(\iota)}=\psi|_{N(\iota)}$.
\item The map $\wt f:V\to T^\perp M$ with $\wt f|_{\gamma(1)}=\tau_\gamma|_{V_{\gamma(1)}}$ is well defined where $T^\perp M$ is the normal bundle of the isometric immersion $f:M\to \wt M$. Moreover, $\wt f$ preserves metrics and connections.
\item  $\wt f^*h_{\wt M}=h$ where $h_{\wt M}$ is the second fundamental form of the isometric immersion $f:M\to \wt M$.
\end{enumerate}
\end{thm}

Although the assumptions of Theorem \ref{thm-CAH-immersion} look rather complicated, they are in fact essentially necessary for the existence of such isometric extensions. For examples, the assumptions (1) , (2) and (3) are the Gauss-Codazzi-Ricci equations for isometric immersions; the assumption \eqref{eq-difference-1} comes from the identity:
$$h_{S,\wt M}(X,Y)=h_{S,M}(X,Y)+h_{M,\wt M}(X,Y),\ \forall\ X,Y\in \Gamma(TS)$$
where $h_{S,M}$ means the second fundamental form for $S\hookrightarrow M$.

When $(\wt M,\wt g)$ is the space form $\mathbb M_K$, we have
$$\wt R(X,Y,Z,W)=K\left(\vv<X,W>\vv<Y,Z>-\vv<X,Z>\vv<Y,Z>\right).$$
So,
$$\left(\tau_{\gamma}^*\wt R\right)(X,Y,Z,W)=K\left(\vv<X,W>\vv<Y,Z>-\vv<X,Z>\vv<Y,Z>\right)$$
which is independent of $\gamma$ since $\tau_\gamma$ is a linear isometry. Therefore, in this case, Theorem \ref{thm-CAH-immersion} has a much simpler form.
\begin{cor}\label{cor-CAH-immersion}
Let $(M^n,g)$ and $(S^r,g_S)$ be two connected Riemannian manifolds, and $\iota:S\to M$ be an isometric immersion with second fundamental $\sigma$ that induces a surjective homomorphism between fundamental groups of $S$ and $M$.  Let $\wt\iota:S\to \mathbb M_K^{n+s}$ be an isometric immersion with second fundamental form $\wt\sigma$. Let $(V^s,\mathfrak{h}, D)$ be a smooth vector bundle on $M$ equipped with the Riemannian metric $\mathfrak{h}$ and the compatible connection $D$, and $h\in \Gamma(\Hom(TM\odot TM, V))$.   Let $\psi:N(\iota)\oplus \iota^*V\to N(\wt\iota)$ be a bundle map  that preserves metrics, connections and satisfies
\begin{equation*}
\psi^*\wt\sigma=\sigma+\iota^*h.
\end{equation*}
Moreover, suppose that for any $X,Y,Z,W\in \Gamma(TM)$ and $\xi,\eta\in \Gamma(V)$,\\
(1)
\begin{equation*}
\begin{split}
R(X,Y,Z,W)=&K\left(\vv<X,W>\vv<Y,Z>-\vv<X,Z>\vv<Y,Z>\right)\\
&+\vv<h(X,W),h(Y,Z)>-\vv<h(X,Z),h(Y,W)>
\end{split}
\end{equation*}
where $R$ is the curvature tensors of $M$;\\
(2)
$$(D_{X}h)(Y,Z)-(D_Yh)(X,Z)=0;$$
(3) $$R^V(\xi,\eta,X,Y)=\vv<A_\xi(Y),A_\eta(X)>-\vv<A_\eta(Y),A_\xi(X)>,$$
 where $R^V$ is the curvature tensor of the vector bundle $V$.

Then, there is an isometric immersion $f:M\to \mathbb M_K^{n+s}$ such that $\wt\iota=f\circ\iota$ and $f_*|_{N(\iota)}=\psi|_{N(\iota)}$.
\end{cor}

Finally, note that the local version of the Cartan-Ambrose-Hicks theorem is Cartan's isometry theorem (see \cite{Ca,CE}), which requires simpler assumptions. So, we also provide the local versions of Theorem \ref{thm-CAH-isometry} and Theorem \ref{thm-CAH-immersion}. We first introduce the following notion.
\begin{defn}
Let $V$ be a smooth vector bundle on $M$. A subset  $\Omega$  of $V$ is said to be star-shaped centered at the zero section if for any $x\in M$, $\Omega_x:=V_x\cap \Omega$ is a star-shaped  subset of $V_x$ centered at the origin.
\end{defn}
\begin{thm}[Local version of Theorem \ref{thm-CAH-isometry}]\label{thm-Cartan-1}
Let $(M^n,g)$ and $(\wt M^n,\wt g)$ be two Riemanian manifolds with $\wt M$ complete.  Let $S^r$ be an embedded submanifold of $M$ and $\wh\Omega$ be a star-shaped open subset of $T^\perp S$ centered at the zero section such that $\exp_S:\wh\Omega\to \Omega\left(:=\exp_S(\wh\Omega)\right)$ is a diffeomorphism. Let $\wt S^r$ be an embedded submanifold of $\wt M$ and $\varphi:S\to \wt S$ be a local isometry. Let $\psi:T^\perp S\to T^\perp \wt S$ be a bundle map along $\varphi$ preserving metrics, second fundamental forms and normal connections. For any $v\in \wh\Omega$, let $\gamma_v(t)=\exp_S(tv)$ and $\wt \gamma_v=\exp_{\wt S}(t\psi(v))$ with $t\in [0,1]$. Suppose that for any $v\in \wh\Omega$ and any plane $\pi\subset T_{\gamma_v(1)}M$ containing $\gamma_v'(1)$, we have
\begin{equation}\label{eq-R}
K(\pi)=K(\tau_v(\pi))
\end{equation}
where $K(\pi)$ means the sectional curvature w.r.t. the plane $\pi$, and
$$\tau_v=P_0^1(\wt\gamma_v)\circ\left(\varphi_*+\psi\right)\circ P_1^0(\gamma_v):T_{\gamma_v(1)}M\to T_{\wt\gamma_v(1)}\wt M.$$
Then, the map $f=\exp_{\wt S}\circ \psi\circ (\exp_S)^{-1}:\Omega\to \wt M$ is the local isometry with $f|_S=\varphi$ and $f_*|_{T^\perp S}=\psi$.
 \end{thm}

\begin{thm}[Local version of Theorem \ref{thm-CAH-immersion}]\label{thm-Cartan-2}
Let $(M^n,g)$ and $(\wt M^{n+s},\tilde g)$ be two Riemannian manifolds with $\wt M$ complete. Let $S^r$ be an embedded submanifold of $M$ with second fundamental form $\sigma$ and $\wh\Omega$ be a star-shaped open subset of $T^\perp S$ centered at the zero section such that $\exp_S:\wh\Omega\to \Omega\left(:=\exp_S(\wh\Omega)\right)$ is a diffeomorphism. Let $(V^s,\mathfrak{h}, D)$ be a smooth vector bundle on $\Omega$ equipped  with the Riemannian metric $\mathfrak{h}$ and  the compatible connection $D$, and $h\in \Gamma(T\Omega\odot T\Omega, V)$. Let $\wt S^r$ be an embedded  submanifold of $\wt M$ with second fundamental forms $\wt \sigma$. Let $\varphi:S\to \wt S$  be a local isometry and $\psi:T^\perp S\oplus V|_S\to T^\perp\wt S$ be a bundle map along $\varphi$ that preserves metrics, connections and satisfies
\begin{equation*}
\psi^*\wt\sigma=\sigma+h|_S.
\end{equation*}
For any $v\in \wh\Omega$, let $\gamma_v(t)=\exp_S(tv)$, $\wt v=\psi(v)$, and $\wt\gamma_v=\dev(\varphi(\gamma(0)), \wt v, \wt h_v)$ where
$$\wt h_{v}(t)=\left(\wh\psi^{-1}\right)^*\left(P_t^0(\gamma_v)h\right), \ \forall t\in [0,1],$$
with $\wh\psi=\varphi_*+\psi$. Moreover, suppose that for any $v\in \wh\Omega$,\\
(1) for any $X,Y,Z,W\in T_{\gamma_v(1)}M$,
$$R(X,Y,Z,W)=(\tau_v^*\wt R)(X,Y,Z,W)+\vv<h(X,W),h(Y,Z)>-\vv<h(X,Z),h(Y,W)>,$$
where $R$ and $\wt R$ are the curvature tensors of $M$ and $\wt M$ respectively;\\
(2)  for any $X,Y,Z\in T_{\gamma_v(1)}M$ and $\xi\in V_{\gamma_v(1)}$,
$$\vv<(D_{X}h)(Y,Z)-(D_Yh)(X,Z),\xi>=(\tau^*_v\wt R)(Z,\xi,X,Y);$$
(3) for any $X,Y\in T_{\gamma_v(1)}M$ and $\xi,\eta\in V_{\gamma_v(1)}$,
$$R^V(\xi,\eta,X,Y)=(\tau_{v}^*\wt R)(\xi,\eta,X,Y)+\vv<A_\xi(Y),A_\eta(X)>-\vv<A_\eta(Y),A_\xi(X)>,$$
 where $R^V$ is the curvature tensor of the vector bundle $V$.
Here
    $$\tau_v=D_0^1(\wt\gamma_v)\circ \wh\psi\circ P_1^0(\gamma):T_{\gamma_v(1)}M\oplus V_{\gamma_v(1)}\to T_{\wt\gamma_v(1)}\wt M.$$
Then, we have the following conclusions.
\begin{enumerate}
\item The map $f(\gamma_v(1))=\wt\gamma_v(1)$ from $\Omega$ to $\wt M$ is an isometric immersion such that $f|_S=\varphi$ and $f_*|_{T^\perp S}=\psi|_{T^\perp S}$.
\item The map $\wt f:V\to T^\perp \Omega$ with $\wt f|_{\gamma_v(1)}=\tau_v|_{V_{\gamma_v(1)}}$ preserves metrics and connections where $T^\perp \Omega$ is the normal bundle of the isometric immersion $f:\Omega\to \wt M$.
\item  $\wt f^*h_{\wt M}=h$ where $h_{\wt M}$ is the second fundamental form of the isometric immersion $f:\Omega\to \wt M$.
\end{enumerate}
\end{thm}

When the target space $(\wt M,\wt g)$ is a space form. We have the following local versions of Corollary \ref{cor-rigidity-isometry} and Corollary \ref{cor-CAH-immersion} which are direct corollaries of Theorem \ref{thm-Cartan-1} and Theorem \ref{thm-Cartan-2}.
\begin{cor}[Local version of Corollary \ref{cor-rigidity-isometry}]
Let $(M^n,g)$ be a Riemanian manifolds $S^r$ be an embedded submanifold of $M$ and $\wh\Omega$ be a star-shaped open subset of $T^\perp S$ centered at the zero section such that $\exp_S:\wh\Omega\to \Omega\left(:=\exp_S(\wh\Omega)\right)$ is a diffeomorphism. Let $\wt S^r$ be an embedded submanifold of $\mathbb M_K^n$ and $\varphi:S\to \wt S$ be a local isometry. Let $\psi:T^\perp S\to T^\perp \wt S$ be a bundle map along $\varphi$ preserving metrics, second fundamental forms and normal connections. For any $v\in \wh\Omega$, let $\gamma_v(t)=\exp_S(tv)$. Suppose that for any $v\in \wh\Omega$ and any plane $\pi$ in $T_{\gamma_v(1)}M$ containing $\gamma_v'(1)$,
$$K(\pi)=K$$
where $K(\pi)$ means the sectional curvature of $g$ w.r.t. the plane $\pi$. Then, the map $$f=\exp_{\wt S}\circ \psi\circ (\exp_S)^{-1}:\Omega\to \mathbb M_K^n$$ is the local isometry with $f|_S=\varphi$ and $f_*|_{T^\perp S}=\psi$.
\end{cor}
\begin{cor}[Local version of Corollary \ref{cor-CAH-immersion}]
Let $(M^n,g)$  be a Riemannian manifold, $S^r$ be an embedded submanifold of $M$ with second fundamental form $\sigma$ and $\wh\Omega$ be a star-shaped open subset of $T^\perp S$ centered at the zero section such that $\exp_S:\wh\Omega\to \Omega:=\exp_S(\wh\Omega)$ is a diffeomorphism. Let $(V^s,\mathfrak{h}, D)$ be a smooth vector bundle on $\Omega$ equipped  with the Riemannian metric $\mathfrak{h}$ and  the compatible connection $D$, and $h\in \Gamma(T\Omega\odot T\Omega, V)$. Let $\wt S^r$ be an embedded  submanifold of $\mathbb M_{K}^{n+s}$ with second fundamental forms $\wt \sigma$. Let $\varphi:S\to \wt S$  be a local isometry and $\psi:T^\perp S\oplus V|_S\to T^\perp\wt S$ be a bundle map along $\varphi$ that preserves metrics, connections and satisfies
\begin{equation*}
\psi^*\tilde\sigma=\sigma+h|_S.
\end{equation*}
Moreover, suppose that for any $X,Y,Z,W\in \Gamma(TM)$ and $\xi,\eta\in \Gamma(V)$,\\
(1)
\begin{equation*}
\begin{split}
R(X,Y,Z,W)=&K\left(\vv<X,W>\vv<Y,Z>-\vv<X,Z>\vv<Y,Z>\right)\\
&+\vv<h(X,W),h(Y,Z)>-\vv<h(X,Z),h(Y,W)>
\end{split}
\end{equation*}
where $R$ is the curvature tensor of $M$;\\
(2)
$$(D_{X}h)(Y,Z)-(D_Yh)(X,Z)=0;$$
(3) $$R^V(\xi,\eta,X,Y)=\vv<A_\xi(Y),A_\eta(X)>-\vv<A_\eta(Y),A_\xi(X)>,$$
 where $R^V$ is the curvature tensor of the vector bundle $V$.

Then, there is an isometric immersion $f:\Omega\to \mathbb M_K^{n+s}$ such that $f|_S=\varphi$ and $f_*|_{T^\perp S}=\psi|_{T^\perp S}$.
\end{cor}
For simplicity, we will adopt the Einstein summation convention for repeated indices and the following conventions of notations throughout the paper:
\begin{enumerate}
\item denote indices in $\{1,2,\cdots,n+s\}$ as $A,B,C,D$;
\item denote indices in $\{1,2,\cdots,n\}$ as $a,b,c,d$;
\item denote indices in $\{1,2,\cdots,r\}$ as $i,j$ ;
\item denote indices in $\{r+1,r+2,\cdots,n\}$ as $\mu,\nu$;
\item denote indices in $\{n+1,n+2,\cdots,n+s\}$ as $\alpha,\beta$;
\item the symbol $'$ means taking derivative with respect to $t$.
\end{enumerate}

The rest of the paper is organized as follows. In section 2, we prove Theorem \ref{thm-CAH-isometry} and Theorem \ref{thm-Cartan-1}. In Section 3, we prove Theorem \ref{thm-CAH-immersion} and Theorem \ref{thm-Cartan-2}. In the appendix, we give a characterization of a smooth map $\iota:S\to M$ inducing a surjective homomorphisms on fundamental groups.


\section{Cartan-Ambrose-Hicks theorem based on submanifolds}
In this section, we will prove Theorem \ref{thm-CAH-isometry}, Theorem \ref{thm-Cartan-1} and Corollary \ref{cor-rigidity-isometry}.

We first derive the equation for the variation field of a family of developments of curves starting from an immersed submanifold.

\begin{lem}\label{lem-eq-J-1}
Let $(M^n,g)$ and $(S^r,g_S)$ be two Riemannian manifolds and $\iota:S\to M$ be an isometric immersion with second fundamental form $\sigma$. Let $\theta(u):I\to S$ be a smooth curve in $S$ and $v(u,t):I\times [0,1]\to TM$ be a smooth map such that $v(u,t)\in T_{\iota\circ\theta(u)}M$ where $I$ is some interval. Let $\gamma_u=\dev(\iota\circ\theta(u),v(u,\cdot))$ and $\Phi(u,t)=\gamma_u(t)$.
Moreover, let $e_1,e_2,\cdots, e_r$ be an orthonormal frame parallel along $\theta$ on $S$ and $e_{r+1}, e_{r+2},\cdots,e_n$ be an orthonormal frame parallel along $\theta$ on the normal bundle $T^\perp S$. Suppose that $v=v_ae_a$, $\theta'=\theta_ie_i$ and
$$\sigma_{ij}^\mu=\vv<\sigma(e_i,e_j),e_\mu>.$$ Let $$E_a(u,t)=P_0^t(\gamma_u)(e_a(u)).$$
Suppose that
\begin{equation*}
\frac{\p \Phi}{\p u}=U_a E_a\mbox{ and }\nabla_\pu E_a=X_{ab}E_b.
\end{equation*}
Then, $X_{ab}=-X_{ba}$, and $U$ and $X$ satisfy the following Cauchy problem:
\[
\begin{cases}
U''_a=R_{bacd}v_bv_cU_d+\p_tv_bX_{ba}+\p_u\p_tv_a\\
X'_{ab}=R_{abcd}v_cU_d\\
U_i(u,0)=\theta_i(u)\\
U_\mu(u,0)=X_{ij}(u,0)=X_{\mu\nu}(u,0)=0\\
X_{i\mu}(u,0)=\sigma_{ij}^\mu\theta_j(u)\\
U'_i(u,0)=\p_uv_i(u,0)-v_\mu(u,0) \sigma_{ij}^{\mu}\theta_j(u)\\
U'_\mu(u,0)=\p_uv_\mu(u,0)+v_i(u,0)\sigma_{ij}^\mu\theta_j(u).
\end{cases}
\]
Here $R_{abcd}=R(E_a,E_b,E_c,E_d)$ with $R$ the curvature tensor of $(M,g)$.
\end{lem}
\begin{proof}
The proof of the ODEs in the Cauchy problem is the same as the proof of Lemma 2.3 in \cite{Yu}. We only need to show the initial data in the Cauchy problem.

Note that
\begin{equation*}
U_a(u,0)e_a=\frac{\p \Phi}{\p u}\bigg|_{t=0}=\theta'(u)=\theta_i(u)e_i.
\end{equation*}
So, $U_i(u,0)=\theta_i(u)$ and $U_\mu(u,0)=0$. Because $e_i$ is parallel along $\theta$ on $S$, we  have
\begin{equation*}
X_{ia}(u,0)e_a=\nabla_{\pu}E_i\bigg|_{t=0}=\nabla_{\theta'}e_i=\nabla^\top_{\theta'}e_i+\sigma(\theta',e_i)=\sigma_{ij}^\mu\theta_j e_\mu.
\end{equation*}
So, $X_{ij}(u,0)=0$ and $X_{i\mu}(u,0)=\sigma_{ij}^\mu\theta_j(u)$. Similarly, because $e_\mu$ is parallel along $\theta$ on $T^\perp S$, we have
\begin{equation*}
X_{\mu a}(u,0)e_a=\nabla_{\pu}E_\mu\bigg|_{t=0}=\nabla_{\theta'}e_\mu=\nabla^\perp_{\theta'}e_\mu-A_{e_\mu}(\theta')=-\sigma_{ij}^\mu\theta_je_i
\end{equation*}
where $A$ is the shape operator of $S$. Thus, $X_{\mu\nu}(u,0)=0$. Finally, note that
\begin{equation*}
\begin{split}
U'_a(u,0)e_a=&\nabla_{\pt}\pu\bigg|_{t=0}=\nabla_{\pu}\pt\bigg|_{t=0}=\nabla_{\pu}(v_aE_a)\bigg|_{t=0}\\
=&\left(\p_uv_a(u,0)+v_b(u,0)X_{ba}(u,0)\right)e_a.
\end{split}
\end{equation*}
So,
\begin{equation*}
U'_i(u,0)=\p_uv_i(u,0)-v_\mu X_{i\mu}(u,0)=\p_uv_i(u,0)-v_\mu(u,0) \sigma_{ij}^{\mu}\theta_j(u)
\end{equation*}
and
\begin{equation*}
U'_{\mu}(u,0)=\p_uv_\mu(u,0)+v_iX_{i\mu}(u,0)=\p_uv_\mu(u,0)+v_i(u,0)\sigma_{ij}^\mu\theta_j(u).
\end{equation*}
This completes the proof of the lemma.
\end{proof}
Because the proof of Theorem \ref{thm-Cartan-1} is simpler, we first prove Theorem \ref{thm-Cartan-1}.
\begin{proof}[Proof of Theorem \ref{thm-Cartan-1}]
Let $\theta(u):(-\e,\e)\to S$ and $v(u):(-\e,\e)\to \wh\Omega$ be smooth curves such that $v(u)\in T_{\theta(u)}^\perp S$ for $u\in (-\e,\e)$. Let $\gamma_u(t)=\exp_{\theta(u)}(tv(u))$ for $t\in [0,1]$ and $\Phi(u,t)=\gamma_u(t)$.
Let $\wt\theta(u)=\varphi(\theta(u))$, $\wt v(u)=\psi(v(u))$, $\wt\gamma_u(t)=\exp_{\wt\theta(u)}(t\wt v(u))$ and $\wt\Phi(u,t)=\wt\gamma_u(t)$.

Let $e_1,e_2,\cdots, e_r$ and $e_{r+1}, e_{r+2},\cdots,e_n$ be  orthonormal frames that are parallel along $\theta$ on $TS$ and $T^\perp S$ respectively. Let $\wt e_i=\varphi_*(e_i)$ for $i=1,2,\cdots, r$. Because $\varphi$ is a local isometry from $S$ to $\wt S$, $\wt e_1,\wt e_2,\cdots, \wt e_r$ is an orthonormal frame parallel along $\wt\theta$ on $\wt S$. Let $\wt e_\mu=\psi(e_\mu)$ for $\mu=r+1,\cdots,n$. Because $\psi$ preserves metrics and connections, $\wt e_{r+1},\cdots,\wt e_n$ is an orthonormal frame parallel along $\wt\theta$ on the normal bundle $T^\perp \wt S$.

Let $v(u)=v_a(u)e_a$ and $\wt v(u)=\wt v_{a}(u)\wt e_a$. It is then clear that $v_a(u)=\wt v_a(u)$ for $a=1,2,\cdots,n$. Let $\theta'=\theta_ie_i$ and $\wt\theta'=\wt\theta_i\wt e_i$. It is also clear that $\theta_i=\wt\theta_i$. Moreover, let $\sigma_{ij}^\mu=\vv<\sigma(e_i,e_j),e_\mu>$ and $\wt\sigma_{ij}^\mu=\vv<\wt\sigma(e_i,e_j),e_\mu>$. Because $\psi$ preserves metrics and second fundamental forms, $\sigma_{ij}^\mu=\wt\sigma_{ij}^\mu$. Here $\sigma$ and $\wt\sigma$ are the second fundamental forms of $S$ and $\wt S$ respectively.

Let $E_a(u,t)=P_0^t(\gamma_u)(e_a)$ and $\wt E_a(u,t)=P_0^t(\wt\gamma_u)(\wt e_a)$. Suppose that
$$\frac{\p\Phi}{\p u}=U_aE_a\mbox{ and }\frac{\p\wt\Phi}{\p u}=\wt U_aE_a.$$
Noting that $v_a=\wt v_a$ are independent of $t$, by Lemma \ref{lem-eq-J-1}, we know that $U$ and $\wt U$ satisfy the following Cauchy problems:
\[
\begin{cases}
U''_a=R_{bacd}v_bv_cU_d\\
U_i(u,0)=\theta_i(u)\\
U_\mu(u,0)=0\\
U'_i(u,0)=\p_uv_i(u,0)-v_\mu(u,0) \sigma_{ij}^{\mu}\theta_j(u)\\
U'_\mu(u,0)=\p_uv_\mu(u,0)+v_i(u,0)\sigma_{ij}^\mu\theta_j(u)
\end{cases}
\]
and
\[
\begin{cases}
\wt U''_a=\wt R_{bacd}\wt v_b\wt v_c\wt U_d\\
\wt U_i(u,0)=\wt\theta_i(u)\\
\wt U_\mu(u,0)=0\\
\wt U'_i(u,0)=\p_u\wt v_i(u,0)-\wt v_\mu(u,0) \wt\sigma_{ij}^{\mu}\wt\theta_j(u)\\
\wt U'_\mu(u,0)=\p_u\wt v_\mu(u,0)+\wt v_i(u,0)\wt\sigma_{ij}^\mu\wt\theta_j(u)
\end{cases}
\]
respectively. Moreover, the curvature assumption \eqref{eq-R} gives us
$$R_{bacd}v_bv_c=\wt R_{bacd}\wt v_b\wt v_c.$$
Thus, $U$ and $\wt U$ satisfy the same Cauchy problem and hence
$$\wt U_a(u,t)=U_a(u,t)$$
by the uniqueness of solutions for Cauchy problems of ODEs.

Note that $f(\gamma_u(t))=\wt\gamma_u(t)$. So,
$$f_*\left(\frac{\p\Phi}{\p u}\right)=\frac{\p\wt\Phi}{\p u}.$$
Because
$$\left\|\frac{\p\Phi}{\p u}\right\|^2=\sum_{a=1}^nU_a^2=\sum_{a=1}^n\wt U_a^2=\left\|\frac{\p\wt\Phi}{\p u}\right\|^2,$$
$f$ is a local isometry. Moreover, since $$f(\theta(u))=\wt\theta(u)=\varphi(\theta(u)),$$
$f|_S=\varphi$, and since
$$f_*(\gamma'_u(0))=\wt\gamma_u'(0)=\psi(\gamma'_u(0)),$$
$f_*|_{T^\perp S}=\psi$. This completes the proof of Theorem \ref{thm-Cartan-1}.
\end{proof}

We next come to prove Theorem \ref{thm-CAH-isometry}.
\begin{proof}[Proof of Theorem \ref{thm-CAH-isometry}]
For each $x\in M$, let $\gamma_0,\gamma_1:[0,1]\to M$ be two smooth curves with $\gamma_0(0)=\iota(p_0)$, $\gamma_1(0)=\iota(p_1)$ for $p_0,p_1\in S$ and $\gamma_0(1)=\gamma_1(1)=x$.  By Proposition \ref{prop-top} in the appendix, there is a smooth curve $\theta:[0,1]\to S$ joining $p_0$ to $p_1$, and a smooth homotopy $\Phi:[0,1]\times[0,1]\to M$ such that
\[
\begin{cases}
\Phi(u,0)=\iota\circ\theta(u)&\forall u\in [0,1]\\
\Phi(0,t)=\gamma_0(t)&\forall t\in [0,1]\\
\Phi(1,t)=\gamma_1(t)&\forall t\in [0,1]\\
\Phi(u,1)=x&\forall u\in[0,1].
\end{cases}
\]
Let $\gamma_u(t)=\Phi(u,t)$ and $v(u,t)=v_{\gamma_u}(t).$
Let $\wt v=\ol\psi(v)$, and
\begin{equation*}
\wt\Phi(u,t)=\wt\gamma_u(t)=\dev(\wt\iota\circ\theta(u),\wt v(u,\cdot))(t).
\end{equation*}

 Moreover, let $e_1,e_2,\cdots,e_r$ be an orthonormal frame parallel along $\theta$ on $S$ and $e_{r+1},e_{r+2},\cdots,e_n$ be an orthonormal frame parallel along $\theta$ on $N(\iota)$. Let $\wt e_a=\ol\psi(e_a)$. Then, $\wt e_{r+1},\wt e_{r+2}, \cdots,\wt e_n$ is an orthonormal frame parallel along $\theta$ on $N(\wt\iota)$ since $\psi$ preserves metrics and normal connections.

Let $v=v_ae_a$, $\theta'=\theta_ie_i$, and
\begin{equation*}
\sigma_{ij}^\mu=\vv<\sigma(e_i,e_j),e_\mu>
\end{equation*}
where $\sigma$ is the second fundamental form of $\iota:S\to M$. Let $\wt v=\wt v_a\wt e_a$ and
\begin{equation*}
\wt\sigma_{ij}^\mu=\vv<\wt \sigma(\wt e_i,\wt e_j),\wt e_\mu>
\end{equation*}
where $\wt\sigma$ is the second fundamental form of $\wt \iota:S\to\wt M$. It is then clear that
\begin{equation}\label{eq-equal}
v_a=\wt v_a,\mbox{ and }\sigma_{ij}^\mu=\wt\sigma_{ij}^\mu
\end{equation}
since $\psi$  preserves metrics and second fundamental forms.

Furthermore, let $E_a(u,t)=P_0^t(\gamma_u)(e_a(u))$ and $\wt E_a(u,t)=P_0^t(\wt\gamma_u)(\wt e_a(u))$. Suppose that
\begin{equation*}
\frac{\p\Phi}{\p u}=U_aE_a \mbox{ and }\frac{\p\wt\Phi}{\p u}=\wt U_a\wt E_a,
\end{equation*}
and
\begin{equation*}
\nabla_{\pu}E_a=X_{ab}E_b\mbox{ and }\wt\nabla_{\pu}\wt E_a=\wt X_{ab}\wt E_b.
\end{equation*}
By the curvature assumption \eqref{eq-R-1}, we have
\begin{equation*}
R_{abcd}=\wt R_{abcd}.
\end{equation*}
So, by Lemma \ref{lem-eq-J-1} and \eqref{eq-equal}, $(U,X)$ and $(\wt U,\wt X)$ satisfy the same Cauchy problem. Therefore
\begin{equation*}
U_a(u,t)=\wt U_a(u,t)\mbox{ and } X_{ab}(u,t)=\wt X_{ab}(u,t).
\end{equation*}
In particular, $\wt U_a(u,1)=U_a(u,1)=0$ and hence
$$\wt \gamma_u(1)=\wt\gamma_0(1),\ \forall\ u\in [0,1].$$
Thus, $f$ is well defined.  By re-parametrization of the curve $\gamma$, it is clear that
$$f(\gamma(t))=\wt \gamma(t),\ \forall t\in [0,1].$$
So
$$U_af_*(E_a)=f_*\left(\frac{\p \Phi}{\p u}\right)=\frac{\p\Phi}{\p u}=\wt U_a\wt E_a=U_a\wt E_a.$$
Thus,
$$f_*(E_a)=\wt E_a,$$
and hence $f$ is a local isometry with $f\circ\iota=\wt\iota$ and $f_*|_{N(\iota)}=\psi$.
\end{proof}
At the end of this section, we come  to prove Corollary \ref{cor-rigidity-isometry}.
\begin{proof}[Proof of Corollary \ref{cor-rigidity-isometry}]
By Theorem \ref{thm-CAH-isometry}, there is a local isometry $f:\mathbb M_K\to \mathbb M_K$ such that $f\circ\iota=\wt\iota$ and $f_*|_{N(\iota)}=\psi$. Since local isometries between complete Riemannian manifolds must be covering maps (see \cite{Hi,CE}) and that $\mathbb M_K$ is simply connected, $f$ must be a diffeomorphism and hence an isometry.
\end{proof}
\section{Cartan-Ambrose-Hicks theorem based on submanifolds for isometric immersions}
In this section, we will prove Theorem \ref{thm-CAH-immersion}. Similar as before, we first come to derive the equation for the variation field of a family of developments of curves w.r.t. symmetric tensors starting from a submanifold.

\begin{lem}\label{lem-eq-J-2}
Let $(S^r,g_S)$, $(M^n,g)$ and $(\wt M^{n+s},\wt g)$ be  Riemannian manifolds with $\wt M$ complete, and $\iota:S\to M$ and $\wt\iota:S\to \wt M$ be isometric immersions with second fundamental forms $\sigma$ and $\wt\sigma$ respectively. Let $(V^s,\mathfrak{h}, D)$ be a smooth vector bundle on $M$ equipped with the Riemannian metric $\mathfrak{h}$ and the compatible connection $D$, and $h\in \Gamma(\Hom(TM\odot TM, V))$. Let $\psi:N(\iota)\oplus \iota^*V\to N(\wt\iota)$ be a bundle map that preserves metrics, connections and satisfies
\begin{equation}\label{eq-sum-2}
\psi^*\wt\sigma=\sigma+\iota^*h.
\end{equation}
Let $\theta(u):I\to S$ be a smooth curves with $I$ some intervals and $\Phi(u,t):I\times [0,1]\to M$ be a smooth map with $\Phi(u,0)=\iota\circ\theta(u)$. Let $\gamma_u(t)=\Phi(u,t)$, $v(u,t)=v_{\gamma_u}(t)$ and $\wt v=\ol\psi(v)$. Let
\begin{equation*}
\wt \Phi(u,t)=\wt \gamma_u(t)=\dev\left(\wt\iota\circ\theta(u), \wt v(u,\cdot),\wt h(u,\cdot)\right)(t)
\end{equation*}
where
 $$\wt h(u,t)=\left(\ol\psi^{-1}\right)^*\left(P_t^0(\gamma_u)(h)\right).$$
Moreover, let
 $e_1,\cdots,e_r$ be an orthonormal frame parallel along $\theta $ on $S$, $e_{r+1},\cdots,e_n$ be an orthonormal frame parallel along $\theta$ on $N(\iota)$, and $e_{n+1},\cdots, e_{n+s}$ be an orthonormal frame parallel along $\theta$ on $V$. Let
 $$\wt e_A=\ol\psi(e_A), E_A(u,t)=P_0^t(\gamma_u)(e_A(u)) \mbox{ and } \wt E_A(u,t)=D_0^t(\wt\gamma_u)(\wt e_A(u)).$$ Suppose that
\begin{equation*}
\theta'(u)=\theta_i(u)e_i(u),\ \sigma(e_i,e_j)=\sigma_{ij}^\mu e_\mu,\ v(u,t)=v_a(u,t) e_a(u)
\end{equation*}
and
\begin{equation*}
h(E_a(u,t),E_b(u,t))=h_{ab}^\alpha(u,t)E_\alpha(u,t).
\end{equation*}
Let
\begin{equation*}
\frac{\p\wt \Phi}{\p u}=\wt U_A\wt E_A\ \mbox{and}\ \wt \nabla_{\pu}\wt E_A=\wt X_{AB}\wt E_B.
\end{equation*}
Then, $\wt X_{AB}=-\wt X_{BA}$, and $(\wt U,\wt X)$ satisfies the following Cauchy problem:
\begin{equation}\label{eq-g-Jacobi-3}
\left\{\begin{array}{rl}\wt U''_a=&2\wt U'_\alpha h_{ab}^\alpha v_b+\wt U_\alpha\p_t(h_{ab}^\alpha v_b)+\wt U_bh_{bc}^\alpha h_{ad}^\alpha v_c v_d+\wt R_{bacA} \wt U_A v_bv_c\\
&+(\p_t v_b)\wt X_{ba}-v_b v_c h_{bc}^\alpha\wt X_{a\alpha }+\p_u\p_tv_a\\
\wt U''_\alpha=&-2\wt U'_a h_{ab}^\alpha v_b-\wt U_{a}\p_t( h_{ab}^\alpha v_b)+
\wt U_\beta h_{ab}^\beta h_{bc}^\alpha  v_a v_c+\wt R_{a\alpha bA}\wt U_A v_av_b\\
&+(\p_t v_a)\wt X_{a\alpha}+v_av_bh_{ab}^\beta\wt X_{\beta\alpha}+\p_u( v_av_bh_{ab}^\alpha)\\
\wt X'_{ab}=&\wt X_{a\alpha}h_{bc}^\alpha v_c-h_{ac}^\alpha v_c\wt X_{b\alpha}+\wt R_{abcA}\wt U_A v_c\\
\wt X'_{a\alpha}=&-\wt X_{ab}h_{bc}^\alpha v_c+ h_{ab}^\beta v_b\wt X_{\beta\alpha}+\wt R_{a\alpha bA}\wt U_Av_b+\p_u(h_{ab}^\alpha v_b)\\
\wt X'_{\alpha\beta}=&\wt X_{a\alpha} h_{ab}^\beta v_b-\wt X_{a\beta}h_{ab}^\alpha v_b+\wt R_{\alpha\beta aA}v_a\wt U_A\\
\wt U_i(u,0)=&\theta_i(u)\\
\wt U_\mu(u,0)=&\wt U_\alpha(u,0)=\wt U'_\alpha(u,0)=\wt X_{ij}(u,0)=\wt X_{\mu\nu}(u,0)=\wt X_{\alpha\beta}(u,0)=0\\
\wt X_{i\mu}(u,0)=&\sigma_{ij}^\mu\theta_j(u)\\
\wt X_{a\alpha}(u,0)=&h_{ai}^\alpha(u,0)\theta_i(u)\\
\wt U'_i(u,0)=&\p_uv_i(u,0)-v_\mu(u,0) \sigma_{ij}^{\mu}\theta_j(u)\\
\wt U'_\mu(u,0)=&\p_uv_\mu(u,0)+v_i(u,0)\sigma_{ij}^\mu\theta_j(u).
\end{array}\right.
\end{equation}
Here $\wt R_{ABCD}=\wt R(\wt E_A,\wt E_B,\wt E_C,\wt E_D)$ with $\wt R$ the curvature tensor of $(\wt M,\wt g)$.
\end{lem}
\begin{proof} The proofs of the ODEs for $\wt U$ and $\wt X$ in \eqref{eq-g-Jacobi-3} are the same as the proof of Theorem 2.2 in \cite{Yu2}. We only need to verify the initial data.

Note that
\begin{equation*}
\wt U_A(u,0) \wt e_A=\frac{\p\wt \Phi}{\p u}\bigg|_{t=0}=\theta_i(u)\wt e_i.
\end{equation*}
So, $\wt U_i(u,0)=\theta_i(u)$ and $\wt U_\mu(u,0)=\wt U_\alpha(u,0)=0$.

By that $e_i$ is parallel along $\theta$ on $S$ and  \eqref{eq-sum-2},
\begin{equation*}
\begin{split}
\wt X_{iA}(u,0)\wt e_A=\wt\nabla_{\pu}\wt E_i\bigg|_{t=0}=\wt\nabla_{\theta'}\wt e_i=\wt\sigma(\wt\theta',\wt e_i)=\sigma_{ij}^\mu \theta_j\wt e_\mu+h^\alpha_{ij}(u,0)\theta_j \wt e_\alpha,
\end{split}
\end{equation*}
we have $\tilde X_{ij}(u,0)=0$, $\wt X_{i\mu}(u,0)=\sigma_{ij}^\mu \theta_j(u)$, and $\wt X_{i\alpha}(u,0)=h^\alpha_{ij}(u,0)\theta_j(u)$.

Moreover, by that $\psi$ preserves metrics and connections, and \eqref{eq-hat-D-1}
\begin{equation*}
\begin{split}
\wt X_{\mu\nu}(u,0)
=&\vv<\wt\nabla_{\pu}\wt E_\mu\bigg|_{t=0},\wt e_\nu>=\vv<\wt\nabla_{\theta'}\wt e_\mu,\wt e_\nu>=\vv<\wt\nabla_{\theta'}^\perp\wt e_\mu,\wt e_\nu>\\
=&\vv<\wh D_{\theta'}e_\mu,e_\nu>=\vv<\nabla^\perp_{\theta'}e_\mu+h(\theta',e_\mu),e_\nu>=0
\end{split}
\end{equation*}
since $e_\mu$ is parallel along $\theta$ w.r.t. the normal connection on $N(\iota)$, where $\wt\nabla^\perp$ means the normal connection on $N(\wt\iota)$. By the same reason,
\begin{equation*}
\begin{split}
\wt X_{\mu\alpha}(u,0)=&\vv<\wt\nabla_{\theta'}\wt e_\mu,\wt e_\alpha>=\vv<\wt\nabla_{\theta'}^\perp\wt e_\mu,\wt e_\alpha>\\
=&\vv<\wh D_{\theta'}e_\mu,e_\alpha>=\vv<\nabla^\perp_{\theta'}e_\mu+h(\theta',e_\mu),e_\alpha>
=h_{i\mu}^\alpha(u,0) \theta_i(u).
\end{split}
\end{equation*}
Similarly, by that $\psi$ preserves metrics and connections, and  that $e_\alpha$ is parallel along $\theta$,
\begin{equation*}
\begin{split}
\wt X_{\alpha\beta}(u,0)
=\vv<\wt\nabla_{\theta'}\wt e_\alpha,\wt e_\beta>
=\vv<\wt\nabla_{\theta'}^\perp\wt e_\alpha,\wt e_\beta>
=\vv<\wh D_{\theta'}e_\alpha,e_\beta>
=\vv<D_{\theta'}e_\alpha-A_{e_\alpha}(\theta'),e_\beta>
=0,
\end{split}
\end{equation*}
where we have used \eqref{eq-hat-D-2}.

Finally, by that
\begin{equation*}
\begin{split}
\wt\nabla_{\pt}\pu\bigg|_{t=0}=&\wt\nabla_{\pt}(\wt U_A\wt E_A)\bigg|_{t=0}\\
=&\wt U'_A(u,0)\wt e_A+\wt U_A(u,0)\wt \nabla_{\pt}\wt E_A\bigg|_{t=0}\\
=&\wt U'_A(u,0)\wt e_A+\wt U_i(u,0)\wt \nabla_{\pt}\wt E_i\bigg|_{t=0}\\
=&\wt U'_A(u,0)\wt e_A+\theta_i(u)h_{ia}^\alpha v_a(u,0)\wt e_\alpha
\end{split}
\end{equation*}
where we have used \eqref{eq-g-dev} in the last equality, and on the other hand,
\begin{equation*}
\begin{split}
\wt\nabla_{\pt}\pu\bigg|_{t=0}=&\wt\nabla_{\pu}\pt\bigg|_{t=0}=\wt\nabla_{\pu}(v_a\wt E_a)\bigg|_{t=0}\\
=&\p_uv_a(u,0)\wt e_a+v_a(u,0)\wt X_{aA}(u,0)\wt e_A.\\
\end{split}
\end{equation*}
So, by comparing the last two equalities, we have
\begin{equation*}
\begin{split}
\wt U'_{i}(u,0)=&\p_uv_i(u,0)+v_a(u,0)\wt X_{ai}(u,0)\\
=&\p_uv_i(u,0)+v_\mu(u,0)\wt X_{\mu i}(u,0)\\
=&\p_uv_i(u,0)-\sigma_{ij}^\mu \theta_j(u)v_\mu(u,0),
\end{split}
\end{equation*}
\begin{equation*}
\begin{split}
\wt U'_{\mu}(u,0)=&\p_uv_\mu(u,0)+v_a(u,0)\wt X_{a\mu}(u,0)\\
=&\p_uv_\mu(u,0)+v_i(u,0)\wt X_{i\mu}(u,0)\\
=&\p_uv_\mu(u,0)+\sigma_{ij}^\mu \theta_j(u)v_i(u,0)\\
\end{split}
\end{equation*}
and
\begin{equation*}
\begin{split}
\wt U'_\alpha(u,0)=&v_a(u,0)\wt X_{a\alpha}(u,0)-\theta_i(u)h_{ia}^\alpha v_a(u,0)=0
\end{split}
\end{equation*}
where we have used $\wt X_{i\mu}(u,0)=\sigma_{ij}^\mu \theta_j(u)$ and $\wt X_{a\alpha}(u,0)=h^\alpha_{ai}\theta_i(u)$ which have been obtained before. This completes the proof of the lemma.
\end{proof}
We are now ready to prove Theorem \ref{thm-CAH-immersion} and Theorem \ref{thm-Cartan-2}. Because the proof of Theorem \ref{thm-Cartan-2} is simpler, we first prove it.
\begin{proof}[Proof of Theorem \ref{thm-Cartan-2}]
Let $\theta(u):(-\e,\e)\to S$ and $v(u):(-\e,\e)\to \wh\Omega$ be smooth curves such that $v(u)\in T_{\theta(u)}^\perp S$ for $u\in (-\e,\e)$. Let $$\gamma_u(t)=\exp_{\theta(u)}(tv(u))\mbox{ for $t\in [0,1]$ and }\Phi(u,t)=\gamma_u(t).$$
Let $e_1,e_2,\cdots, e_r$ be an orthonormal frame parallel along $\theta$ on $S$, $e_{r+1}, e_{r+2},\cdots,e_n$ be an orthonormal frame parallel along $\theta$ on the normal bundle $T^\perp S$, and  $e_{n+1},\cdots,e_{n+s}$ be an orthonormal frame of $V$ parallel along $\theta$.

Let $E_A(u,t)=P_0^t(\gamma_u)(e_A)$, and
$$\frac{\p\Phi}{\p u}=U_ae_a,\ \nabla_{\pu}E_a=X_{ab}E_b\mbox{ and }D_{\pu}E_\alpha=X_{\alpha\beta}E_\beta.$$
Suppose that
\begin{equation*}
\theta'(u)=\theta_i(u)e_i(u),\ \sigma(e_i,e_j)=\sigma_{ij}^\mu e_\mu,\ v(u,t)=v_a(u,t) e_a(u)
\end{equation*}
and
\begin{equation*}
h(E_a(u,t),E_b(u,t))=h_{ab}^\alpha(u,t)E_\alpha(u,t).
\end{equation*}

Then, by the same computation as in \cite[Lemma 2.3]{Yu} and that in Lemma \ref{lem-eq-J-1}, $U$ and $X$ satisfy the following Cauchy problem:
\begin{equation}\label{eq-Jacobi-field}
\begin{split}
\left\{\begin{array}{l}U''_a=R_{cadb}v_cv_dU_b\\
X'_{ab}=R_{abdc}v_dU_c\\
X'_{\alpha\beta}=R^V_{\alpha\beta ab}v_aU_b\\
U_i(u,0)=\theta_i(u)\\
U_\mu(u,0)=X_{ij}(u,0)=X_{\mu\nu}(u,0)=X_{\alpha\beta}(u,0)=0\\
X_{i\mu}(u,0)=\sigma_{ij}^\mu\theta_j(u)\\
U'_i(u,0)=\p_uv_i(u)-v_\mu(u) \sigma_{ij}^{\mu}\theta_j(u)\\
U'_\mu(u,0)=\p_uv_\mu(u)+v_i(u)\sigma_{ij}^\mu\theta_j(u).
\end{array}\right.
\end{split}
\end{equation}
where $R_{abcd}=R(E_a,E_b,E_c,E_d)$ and $R^V_{\alpha\beta ab}=R^V(E_\alpha,E_\beta,E_a,E_b)$.

Moreover, let $\wt\theta=\varphi(\theta)$, $\wt v=\psi(v)$, $\wt e_A=(\varphi_*+\psi)(e_A)$, and
$$\wt\Phi(u,t)=\wt\gamma_u(t)=\dev\left(\wt\theta(u), \wt v(u), \wt h(u,\cdot)\right).$$
where $$\wt h(u,t)=\left(\wh\psi^{-1}\right)^*\left(P_t^0(\gamma_u)(h)\right).$$
Let $\wt E_A(u,t)=D_0^t(\wt\gamma_u)(\wt e_A)$ and suppose that
$$\frac{\p \wt\Phi}{\p u}=\wt U_A\wt E_A\mbox{ and }\wt\nabla_{\pu}\wt E_A=\wt X_{AB}\wt E_B.$$
Then, by Lemma \ref{lem-eq-J-2}, $\wt U$ and $\wt X$ satisfy the Cauchy problem \eqref{eq-g-Jacobi-3}. We then claim that
\begin{equation}\label{eq-solution-1}
\left\{\begin{array}{l}\wt U_a=U_a\\
\wt U_\alpha=0\\
\wt X_{ab}=X_{ab}\\
\wt X_{\alpha\beta}=X_{\alpha\beta}\\
\wt X_{a\alpha}=h_{ab}^\alpha U_b.
\end{array}\right.
\end{equation}
In fact, by the same verification as in the proof of Theorem 1.3 in \cite{Yu2} using the assumptions (1), (2), (3), we know that \eqref{eq-solution-1} satisfies the ODEs in \eqref{eq-g-Jacobi-3}. It is also clear by direct computation that  the initial data in \eqref{eq-g-Jacobi-3} is also satisfied by \eqref{eq-solution-1} by using the initial data in \eqref{eq-Jacobi-field}. Thus, by the uniqueness of solution for Cauchy problems of ODEs, we get the claim.

Furthermore, by that
$$U_af_*(E_a)=f_*\left(\frac{\p \Phi}{\p u}\right)=\frac{\p\Phi}{\p u}=\wt U_a\wt E_a=U_a\wt E_a,$$
we have
$$f_*(E_a)=\wt E_a.$$
Thus $f$ is the local isometry satisfying the properties $f|_S=\varphi$ and $f_*|_{T^\perp S}=\psi$.

Note that $\wt f(E_\alpha)=\wt E_\alpha$ for $\alpha=n+1,\cdots, n+s$. So, $\wt f$ preserves metrics. Moreover, since $f_*(E_a)=\wt E_a$,
\begin{equation*}
\wt f\left(D_\pu E_\alpha\right)=\wt f(X_{\alpha\beta}E_\beta)=X_{\alpha\beta}\wt E_\beta=\left(\wt X_{\alpha A}\wt E_A\right)^\perp=\wt\nabla^\perp_{\pu}\wt E_\alpha=\wt\nabla^\perp_{f_*(\pu)}\wt f(E_\alpha).
\end{equation*}
Here $\wt\nabla^\perp$ means the normal connection of the isometric immersion $f:\Omega\to \wt M$. So $\wt f$ preserves connections.

Finally, by \eqref{eq-solution-1} and $f_*(E_a)=\wt E_a$,
\begin{equation*}
\begin{split}
h_{\wt M}\left(\pu, \wt E_b\right)=\left(\wt \nabla_{\pu}\wt E_b\right)^\perp=\wt X_{b\alpha}\wt E_\alpha=U_ah_{ab}^\alpha \wt E_\alpha=\wt f \left(h\left(\pu,E_b\right)\right).
\end{split}
\end{equation*}
Hence $\wt f^*h_{\wt M}=h$. This completes the proof of theorem.
\end{proof}
Finally, we come to prove Theorem \ref{thm-CAH-immersion}.
\begin{proof}[Proof of Theorem \ref{thm-CAH-immersion}] For each $x\in M$, let $\gamma_0,\gamma_1:[0,1]\to M$ be two smooth curves with $\gamma_0(0)=\iota(p_0),\gamma_1(0)=\iota(p_1)$ for $p_0,p_1\in S$ and $\gamma_0(1)=\gamma_1(1)=x$.  By Proposition \ref{prop-top}, there is a smooth curve $\theta:[0,1]\to S$ joining $p_0$ to $p_1$ and a smooth homotopy $\Phi:[0,1]\times[0,1]\to M$ such that
\[
\begin{cases}
\Phi(u,0)=\iota\circ\theta(u)&\forall u\in [0,1]\\
\Phi(0,t)=\gamma_0(t)&\forall t\in [0,1]\\
\Phi(1,t)=\gamma_1(t)&\forall t\in [0,1]\\
\Phi(u,1)=x&\forall u\in [0,1].
\end{cases}
\]
Let $I=[0,1]$ and the other notations be the same as in Lemma \ref{lem-eq-J-2}.
Suppose that
\begin{equation}
\frac{\p\Phi}{\p u}=U_aE_a,\ \nabla_\pu E_a=X_{ab}E_b\mbox{ and }D_\pu E_\alpha=X_{\alpha\beta}E_\beta.
\end{equation}
Then, by Lemma \ref{lem-eq-J-1}, $U$ and $X$ satisfy the Cauchy problem:
\[
\begin{cases}
U''_a=R_{bacd}v_bv_cU_d+\p_tv_bX_{ba}+\p_u\p_tv_a\\
X'_{ab}=R_{abcd}v_cU_d\\
X'_{\alpha\beta}=R^V_{\alpha\beta ab}v_aU_b\\
U_i(u,0)=\theta_i(u)\\
U_\mu(u,0)=X_{ij}(u,0)=X_{\mu\nu}(u,0)=X_{\alpha\beta}(u,0)=0\\
X_{i\mu}(u,0)=\sigma_{ij}^\mu\theta_j(u)\\
U'_i(u,0)=\p_uv_i(u,0)-v_\mu(u,0) \sigma_{ij}^{\mu}\theta_j(u)\\
U'_\mu(u,0)=\p_uv_\mu(u,0)+v_i(u,0)\sigma_{ij}^\mu\theta_j(u).
\end{cases}
\]
By the same argument as in the proof of Theorem \ref{thm-Cartan-2}, we know that\begin{equation}\label{eq-solution}
\left\{\begin{array}{l}\wt U_a=U_a\\
\wt U_\alpha=0\\
\wt X_{ab}=X_{ab}\\
\wt X_{\alpha\beta}=X_{\alpha\beta}\\
\wt X_{a\alpha}=h_{ab}^\alpha U_b.
\end{array}\right.
\end{equation}
By \eqref{eq-solution}, we know that $\wt U_a(u,1)=U_a(u,1)=0$ and $\wt U_\alpha(u,1)=0$. So $$\frac{\p\wt\Phi}{\p u}(u,1)=0$$ and hence
$\wt \gamma_u(1)=\wt\gamma_0(1)$ for $u\in [0,1]$. Thus, $f$ is well defined. Then, by the same argument as in the proof of Theorem \ref{thm-CAH-isometry}, we know that
$$f_*(E_a)=\wt E_a$$
and $f$ is an isometric immersion with $f\circ\iota=\wt\iota$ and $f_*|_{N(\iota)}=\psi|_{N(\iota)}$. This proves the first conclusion of Theorem \ref{thm-CAH-immersion}.

Next, by that $f_*(E_a)=\wt E_a$, we have
$$\tau_{\gamma_u}|_{V_x}:V_x\to T_{f(x)}^\perp M,\ \left(\tau_{\gamma_u}|_{V_x}\right)(E_\alpha(u,1))=\wt E_\alpha(u,1). $$
Note that
$$\p_uE_\alpha(u,1)=D_{\pu}E_\alpha(u,1)=X_{\alpha\beta}(u,1)E_\beta,$$
and since $\wt X_{\alpha a}(u,1)=-\wt X_{a\alpha}(u,1)=-h_{ab}^\alpha U_b(u,1)=0$ and $X_{\alpha\beta}=\wt X_{\alpha\beta}$,
$$\p_u\wt E_\alpha(u,1)=\wt\nabla_{\pu}\wt E_\alpha(u,1)=\wt X_{\alpha\beta}(u,1)\wt E_\beta+\wt X_{\alpha a}(u,1)\wt E_a=X_{\alpha\beta}(u,1)\wt E_\beta.$$
So,
\begin{equation*}
\begin{split}
\p_u(\tau_{\gamma_u}|_{V_x})(E_\alpha(u,1))=&\p_u(\tau_{\gamma_u}|_{V_x}(E_\alpha(u,1)))-\tau_{\gamma_u}|_{V_x}(\p_u(E_\alpha(u,1)))\\
=&\p_u\wt E_\alpha(u,1)-\tau_{\gamma_u}|_{V_x}(X_{\alpha\beta}(u,1)E_\beta)\\
=&0
\end{split}
\end{equation*}
by the last two equalities. Therefore $\tau_{\gamma_u}|_{V_x}$ is independent of $u\in [0,1]$ and thus $\wt f$ is well-defined. The other properties of $\wt f$ can be verified by the same argument as in the proof of Theorem \ref{thm-Cartan-2}. This completes the proof of the theorem.
\end{proof}
\section{Appendix}
In the appendix, we give a proof of a characterization that a smooth map $\iota:S\to M$ induces a surjective homomorphism on fundamental groups. For a continuous path $c$, we use $[c]$ to denote its homotopic class.

Firstly, we need the following key auxiliary lemma.
\begin{lem}\label{lem-key}
Let $X$ be a topological space, $\theta:[0,1]\to X$ be a continuous path and $\gamma_0,\gamma_1:[0,1]\to X$ be two continuous paths with $\gamma_0(0)=\theta(0)$, $\gamma_1(0)=\theta(1)$ and $\gamma_1(1)=\gamma_0(1)=x$. Then, there is a continuous map $\Phi:[0,1]\times[0,1]\to X$ such that
\begin{equation}\label{eq-homotopy-0}
\left\{\begin{array}{ll}\Phi(u,0)=\theta(u)&t\in [0,1]\\
\Phi(0,t)=\gamma_0(t)&t\in [0,1]\\
\Phi(1,t)=\gamma_1(t)&t\in [0,1]\\
\Phi(u,1)=x&u\in [0,1]\\
\end{array}\right.
\end{equation}
if and only if $[\theta]=[\gamma_0\cdot\gamma_1^{-1}]$.
\end{lem}
\begin{proof}
We first prove the only if part. Let $\Phi:[0,1]\times[0,1]\to X$ be a continuous map such that \eqref{eq-homotopy-0} holds. Then, by the Square Lemma \cite[Lemma 7.17]{Lee}, we have
$$[\gamma_0]=[\Phi(0,\cdot)\cdot\Phi(\cdot,1)]=[\Phi(\cdot,0)\cdot\Phi(1,\cdot)]=[\theta\cdot\gamma_1].$$ So
$$[\theta]=[\gamma_0\cdot\gamma_1^{-1}].$$

We next show the if part. Let $h:[0,1]\times[0,1]\to X$ be a continuous map such that
\begin{equation*}
\begin{split}
h(u,t)=\left\{\begin{array}{ll}\theta(t)& u=0\\
\gamma_0\cdot\gamma_1^{-1}& u=1\\
h(u,0)=\theta(0)\\
h(u,1)=\theta(1).\\
\end{array}\right.
\end{split}
\end{equation*}
Note that
$$h(1,t)=\gamma_0\cdot\gamma_1^{-1}(t)=\left\{\begin{array}{ll}\gamma_0(2t)&0\leq t\leq \frac12\\\gamma_1(2(1-t))&\frac12\leq t\leq 1.
\end{array}\right.$$
Let
\begin{equation*}
\Phi(u,t)=\left\{\begin{array}{ll}h(1,t/2)&u=0,\ 0\leq t\leq 1\\
h(t/u,u)&0<u\leq\frac12,\ 0\leq t\leq u\\
h(1,u+(\frac12-u)(t-u)/(1-u))&0<u\leq \frac12,\ u< t\leq 1\\
h(t/(1-u),u)&\frac12<u<1,\ 0\leq t\leq 1-u\\
h(1,(u-\frac12)(1-t)/u+\frac12)&\frac12<u<1,\ 1-u< t\leq 1\\
h(1,1-t/2)&u=1,\ 0\leq t\leq 1.
\end{array}\right.
\end{equation*}
Then $\Phi:[0,1]\times [0,1]\to X$ is continuous, $$\Phi(0,t)=h(1,t/2)=\gamma_0(t),$$
$$\Phi(1,t)=h(1,1-t/2)=\gamma_1(t)$$
and
$$\Phi(u,0)=h(0,u)=\theta(u).$$ This completes the proof of the lemma.
\end{proof}
Finally, we have the following characterization of a continuous map $f:X\to Y$ inducing a surjective homomorphism on fundamental groups.
\begin{prop}\label{prop-top-0}
Let $f:X\to Y$ be a continuous map between two path connected topological spaces and $p\in X$. Then $f_*:\pi_1(X,p)\to \pi_1(Y,f(p))$ is surjective if and only if for any two continuous paths $\gamma_0,\gamma_1:[0,1]\to Y$ with $\gamma_0(0)=f(p_0)$, $\gamma_1(0)=f(p_1)$ and $\gamma_0(1)=
\gamma_1(1)=x$, there is a continuous path $\theta:[0,1]\to X$ joining $p_0$ to $p_1$, and a continuous map $\Phi:[0,1]\times [0,1]\to Y$ such that
\begin{equation}\label{eq-homotopy}
\left\{\begin{array}{ll}
\Phi(u,0)=f\circ\theta(u)&\forall u\in [0,1]\\
\Phi(0,t)=\gamma_0(t)&\forall t\in [0,1]\\
\Phi(1,t)=\gamma_1(t)&\forall t\in [0,1]\\
\Phi(u,1)=x&\forall u\in [0,1].
\end{array}\right.
\end{equation}
\end{prop}
\begin{proof}
We first prove the if part. Let $c:[0,1]\to Y$ be a continuous path with $c(0)=c(1)=f(p)$ and let $x=c(\frac 12)$. Moreover, let
$$\gamma_0(t)=c\left(\frac{t}{2}\right)\mbox{ and }\gamma_1(t)=c\left(1-\frac t2\right)$$
for $t\in [0,1]$. Then, by assumption, there is a continuous path $\theta:[0,1]\to X$ with $\theta(0)=\theta(1)=p$ and a continuous map $\Phi:[0,1]\times[0,1]\to M$ such that \eqref{eq-homotopy} holds. Then, by Lemma \ref{lem-key}, we know that
$$[f\circ\theta]=[\gamma_0\cdot\gamma_1^{-1}]=[c].$$ Thus, $f_*([\theta])=[c]$ and $f_*$ is surjective.

We next prove the only if part. Let $\gamma_0,\gamma_1:[0,1]\to Y$ be two continuous path with $\gamma_0(0)=f(p_0)$, $\gamma_1(0)=f(p_1)$ and $\gamma_0(1)=
\gamma_1(1)=x$.  Let $\sigma_i:[0,1]\to X$ be a continuous path joining $p$ to $p_i$ for $i=0,1$. By assumption, there is a continuous path $c:[0,1]\to X$ with $c(0)=c(1)=p$ such that $$[f\circ c]=f_*[c]=\left[(f\circ \sigma_0)\cdot\gamma_0\cdot\gamma_1^{-1}\cdot (f\circ\sigma_1)^{-1}\right].$$
Thus
$$[f\circ(\sigma_0^{-1}\cdot c\cdot \sigma_1)]=[(f\circ\sigma_0)^{-1}\cdot f\circ c\cdot f\circ\sigma_1]=[\gamma_0\cdot\gamma_1^{-1}].$$
Let $\theta=\sigma_0^{-1}\cdot c\cdot \sigma_1$. Then, by Lemma \ref{lem-key}, we know that there is a continuous map $\Phi:[0,1]\times[0,1]\to Y$ such that \eqref{eq-homotopy} holds. This completes the proof of the proposition.
\end{proof}
Finally, by standard approximation arguments in \cite[\S2,Chaper 2]{Hir}, we have the following characterization for a smooth map inducing a surjective homomorphism on fundamental groups.
\begin{prop}\label{prop-top}
Let $\iota:S\to M$ be a smooth map between two  connected smooth manifolds and $p\in S$. Then $\iota_*:\pi_1(S,p)\to \pi_1(M,\iota(p))$ is surjective if and only if for any two smooth paths $\gamma_0,\gamma_1:[0,1]\to M$ with $\gamma_0(0)=\iota(p_0)$, $\gamma_1(0)=\iota(p_1)$ and $\gamma_0(1)=
\gamma_1(1)=x$, there is a smooth path $\theta:[0,1]\to S$ joining $p_0$ to $p_1$, and a smooth map $\Phi:[0,1]\times [0,1]\to M$ such that
\begin{equation*}
\left\{\begin{array}{ll}
\Phi(u,0)=\iota\circ\theta(u)&\forall u\in [0,1]\\
\Phi(0,t)=\gamma_0(t)&\forall t\in [0,1]\\
\Phi(1,t)=\gamma_1(t)&\forall t\in [0,1]\\
\Phi(u,1)=x&\forall u\in [0,1].
\end{array}\right.
\end{equation*}
\end{prop}

\end{document}